\documentclass[a4paper,
               headinclude=false,
               headsepline,
               footinclude=false,
               footsepline,
               plainfootsepline,
               abstract=true,
               twoside,
               11pt]{scrartcl}
\newif\ifArbeitsfassung\Arbeitsfassungfalse
\newif\ifVersion\Versionfalse
\newif\ifCopyright\Copyrightfalse
\usepackage{xcolor}
\usepackage{hyperref}
\usepackage{srcltx}
\usepackage{scrlayer-scrpage}
\usepackage{amsmath,amsfonts,amsthm,amssymb}
\usepackage{libertinus}
\usepackage{tikz}
\usetikzlibrary{calc}
\addtokomafont{title}{\rmfamily}
\addtokomafont{section}{\rmfamily}
\addtokomafont{subsection}{\rmfamily}
\addtokomafont{minisec}{\normalfont\rmfamily}

\swapnumbers
\theoremstyle{plain}
\newtheorem{theo}{Theorem}[section]
\newtheorem{lemm}[theo]{Lemma}
\newtheorem{coro}[theo]{Corollary}
\newtheorem{prop}[theo]{Proposition}
\newtheorem{namet}[theo]{\myThmName}

\theoremstyle{definition}
\newtheorem{defi}[theo]{Definition}
\newtheorem{defs}[theo]{Definitions}
\newtheorem{exam}[theo]{Example}
\newtheorem{exas}[theo]{Examples}

\newtheorem{rema}[theo]{Remark}
\newtheorem{rems}[theo]{Remarks}

\newtheorem{named}[theo]{\myThmName}
\newtheorem*{named*}{\myThmName}

\newenvironment{ndef*}[1][\kern-.35em]{\edef\myThmName{#1}\begin{named*}}{\end{named*}}

\newcommand{\coloneqq}{\mathrel{\mathop:}=}
\newcommand{\id}{\mathrm{id}}
\newcommand{\Char}{\operatorname{char}}
\newcommand{\CC}{\mathbb C}
\newcommand{\HH}{\mathbb H}
\newcommand{\RR}{\mathbb R}
\newcommand{\UU}{\mathbb U}
\newcommand{\cB}{\mathcal{B}}
\newcommand{\cL}{\mathcal{L}}
\newcommand{\Aut}[2][]{{\mathrm{Aut}_{#1}}(#2)}
\newcommand{\Gal}[1]{{\mathrm{Gal}}(#1)}
\newcommand{\GL}[2]{\mathrm{GL}(\ifx#1\empty\else#1,\fi#2)}
\newcommand{\PGL}[2]{\mathrm{PGL}(\ifx#1\empty\else#1,\fi#2)}
\newcommand{\SU}[3][]{\mathrm{SU}_{#1}(\ifx#2\empty\else#2,\fi#3)}
\newcommand{\PSU}[3][]{\mathrm{PSU}_{#1}(\ifx#2\empty\else#2,\fi#3)}
\newcommand{\U}[3][]{\mathrm{U}_{#1}(\ifx#2\empty\else#2,\fi#3)}
\newcommand{\EU}[2]{\mathrm{EU}(\ifx#1\empty\else#1,\fi#2)}
\newcommand{\PU}[3][]{\mathrm{PU}_{#1}(\ifx#2\empty\else#2,\fi#3)}
\newcommand{\PEU}[3][]{\mathrm{PEU}_{#1}(\ifx#2\empty\else#2,\fi#3)}
\newcommand{\PgU}[3][]{\mathrm{P\Gamma U}_{#1}(\ifx#2\empty\else#2,\fi#3)}

\newcommand{\PSaU}[2][]{\mathrm{PS\alpha U}_{#1}(#2)}
\newcommand{\trgU}[1][]{\mathrm T\ifx\empty#1\else_{[#1]}\fi}
\newcommand{\GHeis}[3]{\mathrm{GH}(#1,#2,#3)}
\newcommand{\Hsu}[1]{\mathfrak{su}(#1)} 
\newcommand{\Hsos}[1]{\mathfrak{so}^*(#1)} 
\makeatletter
\newcommand{\widebar}[1]{%
  {\mathchoice
    {\vbox
    {\m@th \ialign {##\crcr \noalign {\kern 1\p@ }\kern 1\p@ \hrulefill \crcr
        \noalign {\kern 1\p@ \nointerlineskip }%
        $\hfil \displaystyle {#1}\hfil $\crcr }}}
    {\vbox
    {\m@th \ialign {##\crcr \noalign {\kern 1\p@ }\kern 1\p@ \hrulefill \crcr
        \noalign {\kern 1\p@ \nointerlineskip }%
        $\hfil \textstyle {#1} $\crcr }}}
    {\vbox
    {\m@th \ialign {##\crcr \noalign {\kern 1\p@ }\kern 1\p@ \hrulefill \crcr
        \noalign {\kern 1\p@ \nointerlineskip }%
        $\hfil \scriptstyle {#1}\hfil $\crcr }}}
    {\vbox
    {\m@th \ialign {##\crcr \noalign {\kern 1\p@ }\kern 1\p@ \hrulefill \crcr
        \noalign {\kern 1\p@ \nointerlineskip }%
        $\hfil \scriptscriptstyle {#1}\hfil $\crcr }}}%
    }}
\makeatother
\let\gal\widebar
\newcommand{\rType}[3]{{\text{\normalfont\sffamily#1}}_{#2}^{#3}}
\newcommand{\form}[2]{\langle{#1}|{#2}\rangle}
\newcommand{\proj}[1]{\left\lfloor#1\right\rfloor}
\newcommand{\Gr}[2]{\mathrm{Gr}_{#1}(#2)}
\newcommand{\PG}[2]{\edef\test{#1}\ifx\test\empty{\mathrm{PG}}(#2)\else{\mathrm{PG}}(#1,#2)\fi}
\newcommand{\set}[2]{\left\{{#1}\left|\vphantom{#1}\strut\vphantom{#2}\right.\,
    {#2}\right\}}
\newcommand{\medset}[2]{\left\{\smash{#1}\left|\strut\right.\,
    \strut\smash{#2}\right\}}
\newcommand{\smallset}[2]{\{\smash{#1}\left|\vphantom{}\right.\,\smash{#2}\}}
  \newcount\Hour
   \newcount\Minute
   \def\timenow{\Hour=\time \Minute=\Hour 
        \divide\Hour by 60 \number\Hour:%
        \multiply\Hour by 60 
        \global\advance\Minute by -\Hour%
        \ifnum\Minute<10 0\number\Minute\else\number\Minute\fi}
\let\tOday\today \let\today\empty
\lohead[]{\footnotesize\normalfont\rmfamily \MYtitle}
\lehead[]{\footnotesize\normalfont\rmfamily \MYauthor}
\rehead[]{\footnotesize\normalfont\rmfamily \MYtitle}
\rohead[]{\footnotesize\normalfont\rmfamily \MYauthor}
\lefoot[\pagemark]{\pagemark}
\cefoot[\ifVersion\tiny{\normalfont\ttfamily \jobname.tex,
  \ifArbeitsfassung{\color{red}preliminary }\fi Version of \tOday, \timenow}\fi]%
       {\ifVersion\tiny{\normalfont\ttfamily \jobname.tex,
  \ifArbeitsfassung{\color{red}preliminary }\fi Version of \tOday, \timenow}\fi}
\refoot[\ifCopyright\tiny\copyright~{\normalfont\ttfamily by~\FOOTauthor}\fi]{\ifCopyright\tiny\copyright~{\normalfont\ttfamily by~\FOOTauthor}\fi}
\lofoot[\ifCopyright\tiny\copyright~{\normalfont\ttfamily by~\FOOTauthor}\fi]{\ifCopyright\tiny\copyright~{\normalfont\ttfamily by~\FOOTauthor}\fi}
\cofoot[\ifVersion\tiny{\normalfont\ttfamily \jobname.tex,
  \ifArbeitsfassung{\color{red}preliminary }\fi Version of \tOday, \timenow}\fi]%
       {\ifVersion\tiny{\normalfont\ttfamily \jobname.tex,
  \ifArbeitsfassung{\color{red}preliminary }\fi Version of \tOday, \timenow}\fi}
\rofoot[\pagemark]{\pagemark}
\title{An isomorphism of unitals, and an isomorphism of classical groups}

\author{Markus J. Stroppel}

\makeatletter
  \let\MYauthor\shortauthor 
  \let\MYtitle\shorttitle
  \newcommand\FOOTauthor{the author}%
\makeatother

\begin{document}

\maketitle

\begin{abstract}
\noindent%
An isomorphism between two hermitian unitals is proved, and used to
treat isomorphisms of classical groups that are related to the
isomorphism between certain simple real Lie algebras of types A and D
(and rank~3).
\end{abstract}

\noindent%
In the present paper, we use an isomorphism between two hermitian
unitals to treat isomorphisms of classical groups that are related
to the isomorphism between the simple real Lie algebras of
type~$\rType{A}{3}{\CC,1}$ and~$\rType{D}{3}{\HH}$ (in the notation of
Tits~\cite[pp.\,28, 40]{MR0218489}, Helgason~\cite[X \S\,2.1,
\S\,6.2]{MR514561} denotes the algebras in question by $\Hsu{3,1}$
and~$\Hsos{6}$, respectively).

Our incidence geometric approach complements the algebraic approach
used in~\cite[2.14]{KramerStroppel-JoLT} by a geometric explanation
for the exceptional isomorphism of classical groups.  That algebraic
approach works in much greater generality, including certain
characteristic two cases where the unital over the quaternions
collapses into a line, and cannot be used for our purposes.

\section{Hermitian unitals}

We generalize the notion of finite hermitian unital
(see~\cite[p.\,104]{MR0233275}) to the case of hermitian forms
over infinite (and not necessarily commutative) fields, as follows.

\begin{defs}
  Let\/ $K$ be any (not necessarily commutative) field, and let\/
  $\sigma$ be an anti-automorphism of~$K$, with
  $\sigma^2 = \id \ne \sigma$. %

  If $V$ is a vector space over~$K$, and $h\colon V\times V\to K$ is a
  non-degenerate $\sigma$-hermitian form, we define the set %
  \(%
    U_h \coloneqq \smallset{Kv\in\Gr1V}{v\perp_h v} =
    \smallset{P\in\Gr1V}{P\leqq P^{\perp_h}} %
  \) %
  of \emph{absolute points (with respect to~$h$)}. %
  If $d\coloneqq \dim{V}$ is finite, the hermitian form~$h$
  defines a polarity~$\pi_h$ of the projective space
  $\PG{}{V} \cong \PG{d-1}{K}$ (see~\cite[I,\,\S\,5,
  p.\,9\,ff]{MR0072144}, \cite[II.6, p.\,45\,ff]{MR0333959}). The set
  $U_h$ then consists of all points of $\PG{}{V}$ that are incident
  with their image under that polarity.

  Consider a line $L\in\Gr2V$. If the set $U_h\cap\Gr1L$ of absolute
  points on~$L$ contains more than one point then it is called a
  \emph{block} of~$U_h$. The set of all these blocks is denoted
  by~$\cB_h$. %
  Clearly, any two points of\/~$U_h$ are joined by a unique member
  of\/~$\cB_h$.

  If the form~$h$ has Witt index~$1$, we call $(U_h,\cB_h,\in)$ the
  \emph{hermitian unital with respect to~$h$}. 
\end{defs}

\begin{lemm}\label{blocksAndTangents}
  Assume that\/ $h\colon V\times V\to K$ is a non-degenerate
  $\sigma$-hermitian form of Witt index~$1$. %
  If\/ $h$ is trace-valued %
  then the set of blocks through a given point\/ $P\in U_h$ is
  \[
    \set{U_h\cap\Gr1L}{L\in\Gr2V,P<L\not\leqq P^{\perp_h}} \,.
  \]
\end{lemm}
\begin{proof}
  We write $h\colon V\times V\to K\colon (x,y)\mapsto \form
  xy$. Recall (see~\cite[I, \S\,10, p.\,19]{MR0072144}) that~$h$ is
  trace-valued if, and only if, the set $\smallset{\form vv}{v\in V}$
  is contained in $\smallset{s+s^\sigma}{s\in K}$.
  
  Consider any line $L\in\Gr2V$ through $P \in U_h$.  Then $P = Kv$
  with $v\in V\smallsetminus\{0\}$ such that $\form vv = 0$. %
  If $L\leqq P^{\perp_h}$ then every $w\in L\smallsetminus Kv$ satisfies
  $\form ww \ne 0$ because~$h$ has Witt index~$1$. So~$P$ is the
  unique absolute point in~$L$, and $U_h\cap\Gr1L$ contains no
  block. %
  
  If $L\not\leqq P^{\perp_h}$, we pick any $x\in L\smallsetminus P$; then
  $\form xv\ne0$. Replacing~$x$ by a suitable scalar multiple, we
  achieve $\form xv = -1$. %
  For each $s\in K$, we now have $K(sv+x) \in L$ and
  $\form{sv+x}{sv+x} = \form{sv}{sv}+\form{sv}x+\form{x}{sv}+\form{x}x
  = s\form{v}vs^\sigma+s\form{v}x + \form{x}vs^\sigma + \form{x}x =
  0-s-s^\sigma+\form{x}x$.

  If the form~$h$ is trace-valued, we find~$s$ such that $s+s^\sigma =
  \form{x}x$, and $K(sv+x)$ is a second absolute point on~$L$. So
  $U_h\cap\Gr1L$ is indeed a block in that case. 
\end{proof}

From~\cite[I, \S\,10, p.\,19]{MR0072144} %
we recall that every $\sigma$-hermitian form over a field~$K$ with
$\Char{K}\ne2$ is trace valued.  Also, if $\sigma$ acts non-trivially
on the center of~$K$ (in particular, if~$K$ is commutative) then every
$\sigma$-hermitian form is trace-valued.

\begin{exas}
  Let $C|R$ be a separable quadratic extension of commutative fields,
  and let~$\sigma$ be the generator of $\Gal{C|R}$. %
  Then the form %
  \[ %
  h \colon C^3\times C^3 \to C \colon
  \bigl((x_0,x_1,x_2),(y_0,y_1,y_2)\bigr) \mapsto
  x_0^{}y_2^\sigma+x_1^{}y_1^\sigma+x_2^{}y_0^\sigma  %
  \] %
  is not degenerate, trace-valued, and has Witt index~$1$. %
  If $C$ is finite of order~$e$ then the hermitian unital
  $(U_h,\cB_h,\in)$ is 
  the finite hermitian unital of order~$e$.
\end{exas}

\begin{defs}\label{def:translationUnital}
  Let\/ $\UU \coloneqq (U_h,\cB_h,\in)$ be the hermitian unital with
  respect to a non-degenerate hermitian form $h\colon V\times V\to K$
  of Witt index~$1$, let\/ $X\in U_h$ be a point of\/~$\UU$, and let\/
  $(P,\cL,I)$ be any incidence structure. %
  A map $\eta\colon U_h\to P$ is called an \emph{isomorphism}
  from~$\UU$ onto $(P,\cL,I)$ if\/~$\eta$ is bijective, for every
  block $B\in\cB_h$ there exists a unique block $B'\in\cL$ with
  $B^\eta = \set{X\in P}{(X,B')\in I}$, and the resulting map
  $\beta\colon\cB_h\to\cL\colon B\mapsto B'$ is a bijection. %
  As usual, an \emph{automorphism} of\/~$\UU$ is an isomorphism
  of\/~$\UU$ onto~$\UU$ itself.

  An automorphism of\/~$\UU$ is called a \emph{translation of~$\UU$
    with center~$X$} if it leaves invariant every block through~$X$.
  We write $\trgU[X]$ for the set of all translations of\/~$\UU$ with
  center~$X$.
\end{defs}

If $h\colon V\times V\to K$ is a $\sigma$-hermitian form of Witt
index~$1$, then clearly the group $\PgU{V}h$ of collineations induced
by semi-similitudes acts by automorphisms of the hermitian unital
$(U_h,\cB_h,\in)$. %
See~\ref{centerXiG} and~\ref{centerXiH} below for examples of
translations. 

\begin{theo}\label{PUtwotrsUnital}
  Consider an anti-automorphism~$\sigma$ of a (not necessarily
  commutative) field\/~$K$, with $\sigma^2 = \id \ne \sigma$. %
  Let\/ $h\colon V\times V\to K\colon (v,w)\mapsto \form{v}w$ be a
  non-degenerate $\sigma$-hermitian form of Witt index~$1$. %
  If the form is trace-valued (in particular, if\/ $\Char{K}\ne2$ or
  if\/~$K$ is commutative) and $\dim{V}$ is finite then the group\/
  $\PU{V}h$ 
  acts two-transitively on~$U_h$, and thus transitively both
  on\/~$\cB_h$ and on the set of flags of $(U_h,\cB_h,\in)$. %
\end{theo}
\begin{proof}
  As $h$ has Witt index~$1$, there exists $a\in V\smallsetminus\{0\}$ with
  $\form{a}a = 0$, so $Ka$ lies in~$U_h$. As~$h$ is not
  degenerate, there exists $x\in V$ with $\form{a}x \ne 0$. In
  $L\coloneqq Ka+Kx$ there is a second absolute point~$Kb$,
  see~\ref{blocksAndTangents}. 
  
  Let $P,Q$ be two arbitrary points in~$U_h$. Then there are
  $v,w\in V\smallsetminus\{0\}$ with $\form{v}v=0=\form{w}w$ such that
  $P=Kv$ and $Q=Kw$.  As~$h$ has Witt index~$1$, we have
  $\form{v}w \ne 0$. Replacing~$v$ by a suitable scalar multiple, we
  achieve $\form{v}w = 1$. %
  Now Witt's Theorem (see~\cite[\S\,11, p.\,21]{MR0072144}) asserts
  that there exists $A\in\U{V}h$ with $aA = v$ and $bA = w$. The
  induced collineation $\proj{A}\in \PU{V}h$ then maps the pair
  $(Ka,Kb)$ to $(P,Q)$, and maps the block joining~$Ka$ and~$Kb$ to
  the block joining~$P$ and~$Q$.
\end{proof}

\begin{lemm}\label{lem:translationUnital}
  Let\/ $\UU \coloneqq (U_h,\cB_h,\in)$ be the hermitian unital
  with respect to a non-degenerate $\sigma$-hermitian form
  $h\colon V\times V\to K$ of Witt index~$1$.
  \begin{enumerate}
  \item For each point $X\in U_h$, the set\/~$\trgU[X]$ is a subgroup
    of $\Aut\UU$, and a normal subgroup in the stabilizer of~$X$
    in~$\Aut\UU$.
  \item For each block $B\in\cB_h$ through~$X$, the
    subgroup~$\trgU[X]$ acts transitively on the set\/
    $B\smallsetminus\{X\}$. \\%
    In fact, the intersection $\trgU[X] \cap \PU{V}h$ acts
    transitively on that set.
  \end{enumerate}
\end{lemm}
\begin{proof}
  The set~$\trgU[X]$ is the kernel of the action of the stabilizer
  $\Aut\UU_X$ of~$X$ in~$\Aut\UU$ on the set~$\cB_X$ of all blocks
  through~$X$. So~$\trgU[X]$ is a normal subgroup of~$\Aut\UU_X$.

  Pick $v,w\in V$ such that $X=Kv$ and $B = U_h\cap\Gr1L$, where
  $L=Kv+Kw$. Then $\form vv=0$, and without loss of generality, we may
  assume $\form ww=0$ and $\form vw=1$. An easy computation shows that
  $B=\{Kv\}\cup\smallset{K(pv+w)}{p\in K, p+p^\sigma=0}$. For each
  $p\in K$ with $p+p^\sigma=0$, the linear map $M'$ defined by $vM'=v$
  and $wM'= pv+w$ is an isometry of the restriction of~$h$
  to~$L\times L$. As that restriction is not degenerate, the
  space~$L^\perp$ is a vector space complement to~$L$ in~$V$. We
  extend $M'$ to a linear map~$M$ that acts trivially
  on~$L^\perp$. Then~$M$ belongs to $\U{V}h$, and induces a
  collineation $\proj{M}\in\trgU[X]\cap\PU{V}h$ that maps~$Kw$
  to~$K(pv+w)$. This shows that $\trgU[X]\cap\PU{V}h$ is
  transitive on $B\smallsetminus\{X\}$, as claimed.
\end{proof}

\section{Two hermitian forms, and their unitals}

Let\/ $R$ be a commutative field, and let\/ $C|R$ be a quadratic field
extension.  Then the Galois group $\Gal{C|R}$ has order two, and is
generated by an involution $\sigma\colon x\mapsto\gal{x}$. %
We choose an element $i\in C\smallsetminus\{0\}$ with ${i}^\sigma = -i$. %
(If\/ $\Char{R}=2$ then $j$ lies in~$R$; we will exclude that case
later on.) 

We assume that there is an \emph{anisotropic} $\sigma$-hermitian form
on~$C^2$. Without loss of generality (i.e., up to similitude) we may
assume that this form has Gram matrix $N = \left(
  \begin{smallmatrix}
    1 & 0 \\
    0 & s
  \end{smallmatrix}\right)$.

We consider the quaternion field
\[
  H
\coloneqq H_{C|R}^s = \set{\left(
    \begin{matrix}
      a & x \\
      -s\gal{x} & \gal{a}
    \end{matrix}\right)}{a,x\in C} \,.
\]
Using $w \coloneqq \left(
  \begin{smallmatrix}
    0 & 1 \\
    -s & 0 
  \end{smallmatrix}\right)$ and the embedding $c\mapsto \left(
  \begin{smallmatrix}
    c & 0 \\
    0 & \gal{c}
  \end{smallmatrix}\right)$ of~$C$ into~$H$, we obtain $H 
= C+wC$ with
the multiplication rule 
$(a+wb)(c+wd) = ac-s\gal{b}d+w(\gal{a}d+{b}c)$,
for $a,b,c,d\in C$.

\begin{lemm}\label{def:alpha}%
  The map\/ $\alpha \colon a+wb \mapsto \gal{a}+wb$ (where $a,b\in C$)
  is an involutory anti-automorphism of\/~$H$, the fixed points are
  those in $R+wC$.

  We have $(a+wb)+(a+wb)^\alpha = a+\gal{a}+2wb$ %
  and\/ $(a+wb)(a+wb)^\alpha = a\gal{a}-sb\gal{b}+2w\gal{a}b$.
\end{lemm}
\begin{proof}
  In fact, we have $X^\alpha = i^{-1}X^\kappa i$ for each $X\in H$,
  where~$\kappa\colon a+wb \mapsto \gal{a}-wb$ is the standard
  involution of~$H$.  So~$\alpha$ is the composition of an
  anti-automorphism (namely,~$\kappa$) and an (inner) automorphism
  of~$H$. %
  Straightforward calculations yield the remaining assertions.
\end{proof}

We note that $\alpha$ is the standard involution if $\Char{R}=2$. 

\subsection*{A unital in projective space}
\enlargethispage{7mm}%

\begin{defs}\label{def:hermitianFormG}
  On $C^4$, we consider the $\sigma$-hermitian form %
  \[
    g\colon C^4\times C^4\to C\colon
  \bigl((x_0,x_1,x_2,x_3),(y_0,y_1,y_2,y_3)\bigr) \mapsto
  x_0^{}y_3^\sigma+x_3^{}y_0^\sigma+x_1^{}y_1^\sigma+sx_2^{}y_2^\sigma\,.
  \]
  This form has Witt index~$1$ because the norm form of~$H$ is
  anisotropic.
  
  We assume $\Char{R}\ne2$ (so $i\notin R$), and consider
  $\Xi \coloneqq \set{\xi(u,p)}{u\in C^2,p\in Ri} \subseteq
  \PGL4C$, where
  \[
    \xi\bigl((u_0,u_1),p\bigr) \coloneqq \proj{%
      \begin{matrix}
        1 & u_0 & u_1 & p - \frac12{N(u_0+wu_1)}\\[1ex]
        0 & 1 & 0 & -u_0^\sigma \\[1ex]
        0 & 0 & 1 & -su_1^\sigma \\[1ex]
        0 & 0 & 0 & 1
      \end{matrix}
    } \,.
  \]
  (For any matrix $A\in\GL4C$, we denote by $\proj{A}$ the
  corresponding element in~$\PGL4C$, obtained as the coset modulo scalars.)
\end{defs}

\begin{prop}\label{XiTrsUg}
  \begin{enumerate}
  \item We have
    \[
      U_g = \{C(0,0,0,1)\} \cup
      \smallset{C(1,x_1,x_2,x_3)}{x_3^{}+x_3^\sigma =
        -x_1^{}x_1^\sigma-sx_2^{}x_2^\sigma} \,.
    \]
  \item The set\/ $\Xi$ is a subgroup of\/ $\PSU{C^4}g$. That
    subgroup fixes the point\/ $C(0,0,0,1)$, and acts sharply
    transitively on $U_g \smallsetminus \{C(0,0,0,1)\}$. %
    
    In fact, for $u,v\in C^2$ and $p,q\in Ri$ the product in~$\Xi$
    is obtained as
    \( %
    \xi(u,p)\,\xi(v,q) =
    \xi\left(u+v,p+q+\tfrac12{(vMu^\sigma-uMv^\sigma)}\right)
    \), 
    where $M = \left(
      \begin{smallmatrix}
        1 & 0 \\
        0 & s 
      \end{smallmatrix}\right)$.
  \end{enumerate}
\end{prop}
\begin{proof}
  Consider $x = (x_0,x_1,x_2,x_3) \in C^4\smallsetminus\{(0,0,0,0)\}$ with
  $Cx < x^{\perp_g}$. If $x_0=0$ then
  $0 = x_1^{}x_1^\sigma+sx_2^{}x_2^\sigma = N(x_1+wx_2)$, and
  $Cx=(0,0,0,1)$ because the norm form $N$ is anisotropic.  If
  $x_0\ne0$ then we may assume $x_0=1$, and
  $x_3^{}+x_3^\sigma = -x_1^{}x_1^\sigma-sx_2^{}x_2^\sigma$ follows,
  as claimed.

  It is easy to verify $\Xi \subseteqq \SU{C^4}g$, and that each
  element of~$\Xi$ fixes the point $C(0,0,0,1)$. %

  We note
  \(%
    M(u_0,u_1)^\sigma = M\binom{u_0^\sigma}{u_1^\sigma} =
    \binom{u_0^\sigma}{su_1^\sigma}
  \). %
  Straightforward calculations now yield \\
  $N({u_0+wu_1}) = (u_0,u_1)M(u_0,u_1)^\sigma$, %
  and then 
  $-({u+v})M({u+v})^\sigma + 2({vMu^\sigma-uMv^\sigma}) =
  {-uMu^\sigma-vMv^\sigma-uMv^\sigma}$ leads to \\%
  \[%
    \begin{array}{rl}
   \xi(u,p) \, \xi(v,q)  %
    & =  %
    \proj{%
      \begin{matrix}
        1 & u & p - \frac12{N(u_0+wu_1)}\\[0ex]
        0 & E & -Mu^\sigma \vphantom{\frac M2}\\[0ex]
        0 & 0 & 1 \vphantom{\frac12}
      \end{matrix}
    } %
    \proj{%
      \begin{matrix}
        1 & v & q - \frac12{N(v_0+wv_1)}\\[0ex]
        0 & E & -Mv^\sigma  \vphantom{\frac12} \\[0ex]
        0 & 0 & 1 \vphantom{\frac M2}
      \end{matrix}
                } %
      \\[5ex]
      = &      %
    \proj{%
      \begin{matrix}
        1 & u+v & p+q + \frac12({vMu^\sigma-uMv^\sigma}) - \frac12{N(u_0+v_0+w(u_1+v_1))} \\[0ex]
        0 & E & -M(u+v)^\sigma \vphantom{\frac12} \\[0ex]
        0 & 0 & 1 \vphantom{\frac12}
      \end{matrix}
    } \,,%
    \end{array}
    \] %
    where $E = \left(
    \begin{smallmatrix}
      1 & 0 \\
      0 & 1
    \end{smallmatrix}\right)$.
  As $z \coloneqq vMu^\sigma-uMv^\sigma$ satisfies $z+z^\sigma=0$, we
  obtain
  \( \xi(u,p) \, \xi(v,q) =
  \xi\left(u+v,p+q+\frac12{(vMu^\sigma-uMv^\sigma)}\right) \), as
  claimed.  So~$\Xi$ is closed under multiplication. The inverse of
  $\xi(u,p)$ is $\xi(-u,-p) \in \Xi$.

  Finally, we note that $\xi(u,p)$ maps 
  $C(1,0,0,0)$ to
  $C\left(1,u_0,u_1,p - \frac12{N(u_0+wu_1)}\right)$. This
  shows that~$\Xi$ acts sharply transitively on~$U_g$.
\end{proof}

\begin{rems}\label{centerXiG}
  The set $\set{\xi((0,0),p)}{p\in Ri}$ forms both the center and the
  commutator group of the group\/~$\Xi$.  That commutator group is the
  group $\trgU[C(0,0,0,1)]$ of translations of the unital
  $\UU_g = (U_g,\cB_g,\in)$ with center~$C(0,0,0,1)$.

  For the point $C(1,0,0,0) \in U_g$, we obtain
  \[
    \trgU[C(1,0,0,0)] = \set{\proj{%
        \begin{matrix}
          1 & 0 & 0 & 0 \\
          0 & 1 & 0 & 0 \\
          0 & 0 & 1 & 0 \\
          p & 0 & 0 & 1 \\
        \end{matrix}}}%
    {p\in Ri} \,.
  \]
\end{rems}

\goodbreak%
\subsection*{A unital in the quaternion plane}

\begin{defs}\label{def:hermitianFormH}
  We continue to assume $\Char{R}\ne2$. 
  On $H^3$, we consider the $\alpha$-hermitian form %
  \[
    h\colon H^3\times H^3\colon
    \bigl((X_0,X_1,X_2),(Y_0,Y_1,Y_2)\bigr) \mapsto
    X_0^{}Y_2^\alpha+X_1^{}Y_1^\alpha+X_2^{}Y_0^\alpha \,,
  \]
  here $\alpha$ is the involution introduced in~\ref{def:alpha}. %
  The form~$h$ has Witt index~$1$.
  
  We consider the subset
  $\Psi \coloneqq \set{\psi(X,p)}{X\in H, p\in Ri}$ of the
  group~$\PGL3H$, where
    \[
    \psi(X,p) \coloneqq {\proj{%
        \begin{matrix}
          1 & X & p - \frac12{XX^\alpha}\\[1ex]
          0 & 1 & -X^\alpha \\[1ex]
          0 & 0 & 1
        \end{matrix}
      }} \,.
  \]
  (Again, for $A\in\GL3H$, we denote by $\proj{A}$ the
  corresponding element in~$\PGL3H$, obtained as the coset modulo
  central scalars in this case.)
\end{defs}

\enlargethispage{5mm}%
\begin{prop}\label{XiTrsUh}
  \begin{enumerate}
  \item We have
    \[
      U_h = \{H(0,0,1)\} \cup \smallset{H(1,X,Y)}{Y^{}+Y^\alpha =
        -XX^\alpha} \,.
  \]
  \item The set\/ $\Psi$ is a subgroup of\/ $\PU{H^3}h$. That subgroup
    fixes the point\/ $H(0,0,1)$, and acts sharply transitively on $U_h
    \smallsetminus \{H(0,0,1)\}$.

    The multiplication in~$\Psi$ is given by
    \[
      \psi(X,p)\,\psi(Y,q)
      = \psi\left(X+Y,p+q+\tfrac12(YX^\alpha-XY^\alpha)\right) \,.
    \]
  \end{enumerate}
\end{prop}
\begin{proof}
  The proof is quite analogous to the proof of~\ref{XiTrsUg}.
\end{proof}

\begin{rems}\label{centerXiH}
  The center and the commutator group of the group\/~$\Psi$ both coincide with
  $\set{\psi(0,p)}{p\in Ri}$.  That group is the group
  $\trgU[H(0,0,1)]$ of translations of the unital
  $\UU_h = (U_h,\cB_h,\in)$ with center~$H(0,0,1)$.

  For the point $H(1,0,0) \in U_h$, we obtain
  \[
    \trgU[H(1,0,0)] = \set{\proj{%
        \begin{matrix}
          1 & 0 & 0 \\
          0 & 1 & 0 \\
          p & 0 & 1 \\
        \end{matrix}}}%
    {p\in Ri} \,.
  \]
\end{rems}

\begin{rema}
  The groups $\Xi$ and $\Psi$ are examples of generalized
  Heisenberg groups (cp.~\cite{MR1724629}, \cite{MR2926161},
  \cite{MR3535075}). In fact, they are both isomorphic to
  $\GHeis{R^4}{R}{\beta}$, where~$\beta$ is any non-degenerate
  alternating form on~$R^4$. %
  We give a direct isomorphism explicitly, in~\ref{etaXiIso} below.
\end{rema}

\goodbreak%
\subsection*{An isomorphism of unitals}

\begin{defi}\label{def:etaIsoUnitals}
  For each
  $u = (u_0,u_1)\in C^2$ and each $p\in Ri$, we define the point
  \[
    C\left(1,u_0,u_1,p-\tfrac12{N(u_0+wu_1)}\right)^\eta \coloneqq
    H\left(1,u_0+wu_1,p-\tfrac12{(u_0+u_1w)(u_0+u_1w)^\alpha} \right). %
  \]
  Moreover, we put $C(0,0,0,1)^\eta \coloneqq H(0,0,1)$.

  Thus we obtain a bijection
  \( \eta \colon U_g \to U_h \colon P \mapsto P^\eta \), %
  see~\ref{XiTrsUg} and~\ref{XiTrsUh}.
\end{defi}

\begin{theo}\label{etaXiIso}
  We assume $\Char{R}\ne2$, and use the notation introduced in
  \ref{def:alpha}, \ref{def:hermitianFormG}, \ref{def:hermitianFormH},
  and \ref{def:etaIsoUnitals} above.
  \begin{enumerate}
  \item\label{phiHom}%
    The map\/
    \(
      \varphi \colon \Xi \to \Psi \colon %
      \xi\bigl((u_0,u_1),p\bigr) %
      \mapsto \psi\bigl(u_0+wu_1,p\bigr)
    \)
    is an isomorphism of groups.
  \item\label{equivariantPhiEta}%
    For each $u=(u_0,u_1)\in C^2$, each $p\in Ri$, and each point\/
    $P\in U_g$ we have
    $P^{\eta\,\psi(u_0+wu_1,p)} = P^{\xi(u,p)\,\eta}$;
    here~$\eta\colon U_g\to U_h$ is the map introduced
    in~\ref{def:etaIsoUnitals}.
  \item The map $\eta\colon U_g\to U_h$ induces an isomorphism of
    incidence structures from $(U_g,\cB_g,\in)$ onto
    $(U_h,\cB_h,\in)$.
  \end{enumerate}
\end{theo}

\begin{proof}
  We use the multiplication formulae given in~\ref{XiTrsUg}
  and~\ref{XiTrsUh} to prove assertion~\ref{phiHom}. It suffices to
  verify
  \[
    \begin{array}{rl}
      &\,\,\,
        (v_0+wv_1)(u_0+wu_1)^\alpha - (u_0+wu_1)(v_0+wv_1)^\alpha \\
      &=
          (v_0+w{v_1})(\gal{u_0}+w{u_1}) -
          (u_0+w{u_1})(\gal{v_0}+w{v_1})
      \\
      &= v_0\gal{u_0}-u_0\gal{v_0}
          + w^2(\gal{v_1}\,u_1-\gal{u_1}\,v_1) \\
      &= vMu^\sigma-uMv^\sigma;
    \end{array}
  \]
  here we use $wc=\gal{c}w$ (for $c\in C$) and $w^2=-s$.

  \goodbreak%
  Assertion~\ref{equivariantPhiEta} is easily checked. %
  As any two points in a hermitian unital are joined by a unique
  block, it remains to verify that $B^\eta \in \cB_h$ holds for each
  block $B\in\cB_g$. Using transitivity of $\Xi$ on
  $U_g\smallsetminus\{C(0,0,0,1)\}$, we see that it suffices to consider
  blocks through $C(0,0,0,1)$, and blocks through $C(1,0,0,0)$.

  Any block through $C(0,0,0,1)$ is of the form $B = U_g\cap L$, where
  $L = C(0,0,0,1) + C(1,u_0,u_1,u_2)$. %
  We may assume $C(1,u_0,u_1,u_2) \in U_g$. %
  Then $u_2 = p-\frac12{N(u_0+wu_1)}$ holds for some $p\in Ri$.   %
  So the block in question is %
  \[
    B = \set{C(1,u_0,u_1,p-\tfrac12{N(u_0+wu_1)})}{p\in Ri} \cup
    \{C(0,0,0,1)\},
  \]
  and its image
  \[
    B^\eta =
    \set{H(1,u_0+wu_1,p-\tfrac12{(u_0+wu_1)(u_0+wu_1)^\alpha})}{p\in Ri} %
    \cup \{H(0,0,1)\}
  \]
  belongs to~$\cB_h$.

  \goodbreak%
  Now consider a block $B$ through~$C(1,0,0,0)$. %
  There exist $u=(u_0,u_1)\in C^2$ and $x\in Ri$ such that
  $C\left(1,u_0,u_1,x-\frac12 N(u_0+wu_1)\right) \in
  B\smallsetminus\{C(1,0,0,0)\}$.  We abbreviate
  $n \coloneqq N({u_0+wu_1})$. %
  Every point in $B\smallsetminus\{C(1,0,0,0)\}$ is of the form
  $P_a \coloneqq C\left(1,au_0,au_1,a(x-\frac n2)\right)$, where
  $a = a_0+a_1i \in C$ (with $a_0,a_1\in R$) satisfies %
  \[
    a\gal{a}n+2a_1ix-a_0n = 0 \,. %
    \eqno{(*)}%
  \]
  So $P_a = C\left(1,au_0,au_1,y_a-a\gal{a}\frac n2\right)$, with
  $y_a \coloneqq a(x-\frac n2) + a\gal{a}\frac n2$. %
  Note that $y_a \in Ri$.

  We abbreviate $Z \coloneqq u_0+wu_1$, so
  $C(1,u_0,u_1,x-\frac{n}2)^\eta = H(1,Z,x-\frac12ZZ^\alpha)$. %
  For each $a\in C$ satisfying condition~$(*)$ from above, we 
  obtain
  \[
    \begin{array}{rcl}
      C\left(1,au_0,au_1,a({x-\frac{n}2})\right)^\eta
      &=& C\left(1,u_0a,u_1a,y_a-a\gal{a}\frac{n}2\right)^\eta \\[1ex]
      &=& H\left(1,Za,y_a-\frac12Za(Za)^\alpha\right) \\[1ex]
      &=& H\left(1,Za,y_a-\frac12a\gal{a}ZZ^\alpha\right) \,.
    \end{array}
  \]
  Each one of those points is contained in $U_h = U_g^\eta$.  In order
  to see that it is actually contained in the block
  $\left(H(1,0,0)+H(1,Z,x-\frac12ZZ^\alpha)\right) \cap U_h$, it
  remains to check that there exists $Y\in H$ such that
  $Y(Z,x-\frac12YY^\alpha) = (Za,y_a-\frac12a\gal{a}ZZ^\alpha)$. The
  entry on the left yields $Y=ZaZ^{-1}$.  Using
  $ZZ^\alpha = ({u_0+wu_1})({\gal{u_0}+wu_1}) =
  u_0\gal{u_0}-su_1\gal{u_1}+2w\gal{u_0}u_1$ and
  $ZaZ^{-1} = (u_0a+wu_1a)({\gal{u_0}-wu_1})\frac1n =
  a_0+a_1i\,\gal{ZZ^\alpha}\frac1n$, we compute %
  \[ %
    \begin{array}{rcl}
      ZaZ^{-1}(x-\frac12ZZ^\alpha) %
      &=& ZaZ^{-1}x-\frac12ZaZ^\alpha \\
      &=&  a_0x+a_1i\,\gal{ZZ^\alpha}x\frac1n-\frac12a_0ZZ^\alpha-\frac12a_1in \\[1ex]%
      &=& a_0(x-\frac12ZZ^\alpha)+a_1i(\gal{ZZ^\alpha}x\frac1n-\frac12n) \\[1ex]%
      &=&
          a_0(x-\frac12ZZ^\alpha)+a_1i(x{ZZ^\alpha}\frac1n-\frac12n) \,; %
    \end{array}
  \]
  we have used $i\,\gal{F}i^{-1} = F^\alpha$ and $x\in Ri$. %
  On the other hand, we find %
  \[
    \begin{array}{rcl}
      y_a-\frac12a\gal{a}ZZ^\alpha %
      &=& a(x-\frac n2) + a\gal{a}\frac n2 - \frac12a\gal{a}ZZ^\alpha
      \\[1ex]%
      &=& ax-\frac12an+\frac12a_0n-a_1ix-\frac12(a_0-2a_1ix\frac1n)ZZ^\alpha
      \\[1ex]%
      &=& a_0x-a_1i\frac{n}2-\frac12(a_0-2a_1ix\frac1n)ZZ^\alpha \,,
    \end{array}
  \]
  and this equals $ZaZ^{-1}(x-\frac12ZZ^\alpha)$, as required. 

  So we have established that $B^\eta$ is contained in some block~$B'$
  of $\cB_h$, for each $B\in\cB_g$. %
  It remains to show that $B^\eta$ fills all of~$B'$. %
  To this end, we use the fact that the group
  $\trgU \coloneqq \trgU[H(1,0,0)] = \set{\proj{%
      \begin{smallmatrix}
        1 & 0 & 0 \\
        0 & 1 & 0 \\
        p & 0 & 1
      \end{smallmatrix}}}%
  {p\in Ri}$ %
  of translations with center $H(1,0,0)$ acts transitively on
  $D\smallsetminus\{H(1,0,0)\}$, for each block $D\in\cB_h$
  through~$H(1,0,0)$, see~\ref{lem:translationUnital}
  and~\ref{centerXiH}. %
  In particular, we obtain that the block
  $B' = \left(H(1,0,0)+H(1,Z,x-\frac12ZZ^\alpha)\right) \cap U_h$
  equals the set
  $\{H(1,0,0)\} \cup
  \medset{H(1+(x-\frac12ZZ^\alpha)p,Z,x-\frac12ZZ^\alpha}{p\in Ri}$. %
  So it suffices to show that for each $p\in Ri$ there exists $a\in C$
  satisfying condition~$(*)$ and such that
  \[
    H\bigl(1+(x-\tfrac12ZZ^\alpha)p,Z,x-\tfrac12ZZ^\alpha\bigr) =
    H\bigl(1,Za,ZaZ^{-1}(x-\tfrac12ZZ^\alpha)\bigr);
  \]
  the description on the right hand side then yields that the point in
  question lies in~$B^\eta$.

  We need to find $a\in C$ with
  $1+(x-\frac12ZZ^\alpha)p = (ZaZ^{-1})^{-1} = Za^{-1}Z^{-1}$. %
  We write $b \coloneqq a^{-1}$ as $b=b_0+b_1i$ with $b_0,b_1\in R$,
  and compare $1+(x-\frac12ZZ^\alpha)p = (1+xp) -\frac12ZZ^\alpha p$
  with $ZbZ^{-1} = b_0+b_1i\,\gal{ZZ^\alpha}\frac1n$. %
  Since $1+xp$ lies in~$R$ and $\frac12ZZ^\alpha p\in Ri+wC$, we
  obtain $1+xp = b_0$ and
  $-\frac12ZZ^\alpha p = b_1i\,\gal{ZZ^\alpha}\frac1n$, so
  $b_1i = -\frac12pn$, and $b = 1+xp-\frac12pn$.

  Condition~$(*)$ for~$a$ means $n-2b_1ix-b_0n = 0$, and is easily
  verified.  
\end{proof}

\section{Groups of translations, and an isomorphism of
  groups}

\begin{defi}
  Let $(P,\cL,I)$ be an incidence structure such that through any two
  points in~$P$ there is at most one line in~$\cL$ incident with both
  of those points. %
  An \emph{O'Nan configuration} in $(P,\cL,I)$ consists of~$4$ lines
  meeting in~$6$ points (see Fig.~\ref{fig:oNan-fromTranslation}
  below). In particular, any two of those four lines have a (unique)
  point in common.
\end{defi}

These configurations are named in honor of Michael O'Nan, who used
the finite case of the following result~\ref{noONan} in his study of
the automorphisms of finite hermitian unitals, see~\cite{MR0295934}.
In the (axiomatic) context of projective spaces, O'Nan configurations
are called Veblen-Young figures.

The proof of the following result is taken from~\cite[2.2]{MR2795696}.

\begin{prop}\label{noONan}
  Let\/ $V$ be a vector space over a commutative field\/~$F$, and
  assume that there is a non-trivial involutory automorphism~$\sigma$
  of~$F$.  Let\/ $h\colon V\times V\to F\colon (u,v)\mapsto\form uv$
  be a non-degenerate $\sigma$-hermitian form of Witt index~$1$. %
  Then the hermitian unital\/ $\UU=(U_\sigma,\cB_\sigma,\in)$ does
  not contain any O'Nan configurations.
\end{prop}
\begin{proof}
  Consider an O'Nan configuration in the projective space $\PG{}{V}$.
  Then the six points of the configuration are contained in the
  projective plane spanned by any two of the lines inside~$\PG{}{V}$.

  Therefore, there are linearly independent vectors $b_0$, $b_1$,
  $b_2$ in~$V$ such that the six points of the configuration are
  $Fb_0$, $Fb_1$, $F(b_0+b_1)$, $Fb_2$, $F(b_0+b_2)$ and $F(b_1-b_2)$,
  respectively.  If these points belong to~$U_h$ then
  $\form{b_n}{b_n} = 0$ and $\form{b_n}{b_m} = -\form{b_m}{b_n}$ holds
  for all $m < n < 3$.  The matrix $( \form{b_m}{b_n} )_{m,n<3}$ has
  determinant~$0$ (here we use that $F$ is commutative). Hence~$f$ is
  degenerate, and the restriction of~$h$ to $Fb_0+Fb_1+Fb_2$ has Witt
  index at least~$2$. But then the Witt index of~$h$ is greater
  than~$1$, contradicting our assumption.
\end{proof}

\begin{rema}\label{KestenbandCounterexample}
  Kestenband~\cite{MR1187635} claims that~\ref{noONan} holds even for
  hermitian unitals over skew fields. This claim is false. For
  instance, consider the quaternion field
  $\HH \coloneqq H_{\CC|\RR}^1 = \CC+j\CC$ over the real number
  field~$\RR$, constructed from $\CC=\RR+\RR i$ with $j^2=-1$, the
  standard involution~$\kappa\colon x\mapsto\gal{x}$, and the
  hermitian form given by
  \[
    \form{(u_0,u_1,u_2)}{(v_0,v_1,v_2)} = %
    u_0i\gal{v_1}+u_0j\gal{v_2} %
    -u_1i\gal{v_0}-u_1ji\gal{v_2} %
    -u_2j\gal{v_0}+u_2ji\gal{v_1} \,.
  \]
  That form is not degenerate, and has Witt index~$1$.  However, the
  corresponding hermitian unital contains the O'Nan configuration with
  the points $\HH(1,0,0)$, $\HH(0,1,0)$, $\HH(0,0,1)$, $\HH(1,1,0)$,
  $\HH(1,0,1)$, and $\HH(0,1,-1)$.
\end{rema}

\begin{prop}
  Let\/ $\UU = (U_h,\cB_h,\in)$ be a hermitian unital, and
  let\/~$X$ be any point in~$U_h$.  If\/~$\UU$ contains no O'Nan
  configurations then the translation group~$\trgU[X]$ acts sharply
  transitively on $B\smallsetminus\{X\}$, for each block~$B$ through~$X$.
\end{prop}
\begin{proof}
  We already know from~\ref{lem:translationUnital} that~$\trgU[X]$ is
  transitive on $B\smallsetminus\{X\}$.   \begin{figure}[h!]
    \centering
    \begin{tikzpicture}[%
    scale=1, every node/.append style={circle, fill=black,
      inner sep=.5pt, %
      minimum size=3pt}]%
    \node[coordinate] (Y) at (0,0) {} ;%
    \node[coordinate] (X) at (0,2) {} ;%
    \node[coordinate] (W) at (5,.5) {} ;%
    \node[coordinate] (Z) at ($(Y)!.7!(W)$) {} ;%
    \node[coordinate] (Zt) at ($(X)!.7!(Z)$) {} ;%
    \node[coordinate] (Zx) at ($(X)!1.2!(Z)$) {} ;%
    \node[coordinate] (Wx) at ($(X)!1.2!(W)$) {} ;%
    \node[coordinate] (Wt) at (intersection of X--W and Y--Zt) {} ;%
    \node[coordinate] (Dx) at ($(Y)!1.2!(W)$) {} ;%
    \node[coordinate] (Dtx) at ($(Y)!1.4!(Wt)$) {} ;%

    \draw[dashed] (X) -- (Y) ;%
    \draw (X) -- (Zx) ;%
    \draw (X) -- (Wx) ;%
    \draw (Y) -- (Dx) ;%
    \draw (Y) -- (Dtx) ;%

    \node at (X) {} ;%
    \node at (Y) {} ;%
    \node at (Z) {} ;%
    \node at (W) {} ;%
    \node at (Zt) {} ;%
    \node at (Wt) {} ;%

    \node[fill=none,anchor=east] at ($(X)!0.5!(Y)$) {$B$} ;%
    \node[fill=none,anchor=east] at (X) {$X$} ;%
    \node[fill=none,anchor=east] at (Y) {$Y=Y^\tau$} ;%
    \node[fill=none,anchor=north] at (Z) {$Z$} ;%
    \node[fill=none,anchor=north] at (W) {$W$} ;%
    \node[fill=none,anchor=east] at (Zt) {$Z^\tau\quad$} ;%
    \node[fill=none,anchor=south west] at (Wt) {$W^\tau$} ;%
    \node[fill=none,anchor=north west] at (Zx) {$B_Z$} ;%
    \node[fill=none,anchor=north west] at (Wx) {$B_W$} ;%
    \node[fill=none,anchor=west] at (Dx) {$D$} ;%
    \node[fill=none,anchor=west] at (Dtx) {$D^\tau$} ;%
    \end{tikzpicture}
    \caption[Constructing an O'Nan configuration.]{Constructing an
      O'Nan configuration from a translation with a fixed point.}
    \label{fig:oNan-fromTranslation}
  \end{figure}
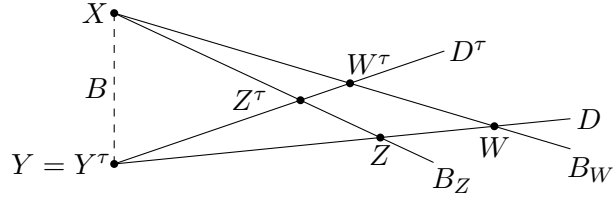
  \\
  If the action is not sharply
  transitive then there exists $\tau\in\trgU[X]\smallsetminus\{\id\}$ such
  that $\tau$ fixes some point $Y\in B\smallsetminus\{X\}$. Let~$Z$ be any
  point in~$U_h\smallsetminus B$ with $Z^\tau\ne Z$. As~$\tau$ is a translation of the
  unital~$\UU$, the block~$B_Z$ joining $X$ and~$Z$ is invariant
  under~$\tau$, and contains~$Z^\tau$. In the block~$D$ joining~$Y$
  and~$Z$, choose a third point~$W$.  Then $W^\tau$ lies in the
  intersection of $D^\tau$ and the block~$B_W$ joining~$X$ and~$W$. %
  So the six points $X$, $Y$, $Z$, $Z^\tau$, $W$, $W^\tau$ and the
  four blocks~$D$, $D^\tau$, $B_Z$, $B_W$ form an O'Nan configuration
  in the unital, contradicting our assumption. See
  Fig.~\ref{fig:oNan-fromTranslation}.
\end{proof}

\begin{coro}
  Let\/ $h\colon V\times V\to K$ be a non-degenerate $\sigma$-hermitian
  form of Witt index~$1$. If the corresponding hermitian unital\/ 
  $\UU$ has no O'Nan configurations then $\trgU[X] =
  \trgU[X]\cap\PU{V}h$ holds for each point of the unital.   %
  \qed
\end{coro}

\goodbreak%
\begin{coro}
  Let\/ $h\colon V\times V\to K$ be a non-degenerate
  $\sigma$-hermitian form of Witt index~$1$. If the corresponding
  hermitian unital\/~$\UU$ has no O'Nan configurations then every
  translation of the unital\/~$\UU$ is induced by a transvection of the
  projective space~$\PG{}{V}$; in fact, each translation with center
  $X=Kv$ is induced by a transvection $\tau_{\lambda,v} \in \U{V}h$
  with $\ker\lambda = v^{\perp_h}$. %
  \qed
\end{coro}

Explicitly, we obtain for the two unitals considered here: %
the commutator groups
$\Xi' = \set{\xi((0,0),p)}{p\in Ri}$ and
$\Psi' = \set{\psi(0,p)}{p\in Ri}$ of~\ref{centerXiG}
and~\ref{centerXiH} are full translation groups, with centers
$C(0,0,0,1)$ and $H(0,0,1)$, respectively.

\begin{theo}\label{PSU4C1isoPaU3H}
  The groups $\PEU{C^4}g$ and $\PEU{H^3}h$ are isomorphic. 
\end{theo}
\begin{proof}
  Recall that the groups $\EU{C^4}g$ and $\EU{H^3}h$, respectively,
  are generated by all unitary transvections; those transvections
  induce the translations of the unital. %
  Conjugation by the isomorphism $\eta\colon \UU_g\to\UU_h$ maps
  $\Aut{\UU_g}$ onto~$\Aut{\UU_h}$, and maps the group $\trgU[X]$
  to~$\trgU[X^\eta]$, for each point $X\in U_g$. %
  So conjugation by~$\eta$ induces an isomorphism from $\PEU{C^4}g$
  onto~$\PEU{H^3}h$.
\end{proof}

\begin{exam}
  We take the field~$\CC$ of complex numbers for~$C$, with the
  standard involution $\sigma \colon c\mapsto \gal{c}$ generating
  $\Gal{\CC|\RR}$, and the field $\HH = H_{\CC|\RR}^1 = \CC+j\CC$ of
  Hamilton's quaternions. %
  The involution~$\alpha$ from~\ref{def:alpha} represents the unique
  class of involutory anti-automorphisms of~$\HH$ apart from the
  standard involution~$\kappa$. %
  For the forms~$g$ and~$h$ introduced
  in~\ref{def:hermitianFormG} and~\ref{def:hermitianFormH},
  respectively, we obtain the groups
  $\PEU{\CC^4}{g} \cong \PSU[4]{\CC}1$ and
  $\PEU{\HH^3}h \cong \PSaU[3]{\HH}$ %
  (in the notation of~\cite[94.33]{MR1384300}, %
  in Tits~\cite[pp.\,28, 40]{MR0218489}, these occur as the groups of
  type~$\rType{A}{3}{\CC,1}$ and~$\rType{D}{3}{\HH}$, %
  Helgason~\cite[X \S\,2.1, \S\,6.2]{MR514561} denotes the
  corresponding algebras by $\Hsu{3,1}$ and~$\Hsos{6}$, respectively).
\end{exam}

\begin{rems}
  For the commutative field~$C$, one knows that $\EU{C^4}g =
  \SU{C^4}g$, so $\PEU{C^4}g = \PSU{C^4}g$. 
  
  Also, it is known that the groups $\PEU{C^4}g$ and $\PEU{H^3}h$ are simple:
  see~\cite[II\,\S\,4]{MR0072144} for a general result,
  cf.~\cite[10.20]{MR1189139} or~\cite[11.26]{MR1859189} for the case
  of a commutative ground field.
  As we restrict our investigation to cases where the characteristic
  is different from two, all the forms in question are trace valued
  forms.
\end{rems}

\begin{rems}
  As the field~$C$ is commutative, the involution $\sigma$ of~$C$ is
  an involution of the second kind (in the sense of
  Dieudonn\'e~\cite[\S\,10, p.\,19]{MR0072144}).  According
  to~\cite[5.6c]{MR2241352}, every reflection in the group $\PU{C^4}g$
  is thus admissible, and we obtain
  $\Aut{U_g,\cB_g,\in} = \PgU{C^4}g$. From our result~\ref{etaXiIso}
  we then also infer $\Aut{U_h,\cB_h,\in} \cong \PgU{C^4}g$.
\end{rems}

\begin{rema}
  Let $F$ be a commutative field, and let $Q$ be a quaternion algebra
  over~$F$. Then $Q$ is a central simple $F$-algebra (cp.~\cite[4.5,
  Lemma\,3, p.\,232]{MR1009787}, and every $F$-linear automorphism is
  inner (by the Skolem-Noether Theorem, see~\cite[p.\,222]{MR1009787},
  or see~\cite[Theorem\,2, p.\,67]{MR0101253} for a direct proof). %
  It then follows that every $F$-linear anti-automorphism~$\beta$ is
  the product of the standard involution and some inner automorphism,
  say $x\mapsto i^{-1}xi$ with $i\in F\setminus\{0\}$, so
  $x^\beta = i^{-1}\gal{x}i$. %
  We obtain that~$\beta$ is an involution precisely if $i^2\in F$,
  i.e., if either $i\in F$ or $\gal{i}=-i$.

  If $i\in F$ then $\beta$ is the standard involution.  If
  $i\notin F$, we form the quadratic extension $C = F+Fi$.  The
  restriction~$\sigma$ of the standard involution of~$Q$ then is the
  generator of~$\Gal{C|F}$, and~$\beta$ is obtained as
  in~\ref{def:alpha}.
\end{rema}

\goodbreak%

\begin{thebibliography}{10}
\providecommand{\href}[2]{#2}
\providecommand{\eprint}[1]{\href{http://arxiv.org/abs/#1}{#1}}
\providecommand{\url}[1]{\href{#1}{#1}}
\providecommand{\urlprefix}{}
\let\oldunderscore_
\catcode`\_=13
\providecommand{\doi}[1]{\href{http://dx.doi.org/#1}{{\def_{\_}\normalfont\ttfamily doi:#1}}\let_\oldunderscore}
\providecommand{\MR}[1]{\relax\ifhmode\unskip\space\fi \MRnumberextract#1 \,}
\def\MRnumberextract#1 #2\,{\MRhref{#1}{#2}}%
\providecommand{\MRhref}[2]{%
  \href{https://mathscinet.ams.org/mathscinet-getitem?mr=#1}{MR\,#1 #2}}
\providecommand{\ZBL}[1]{\relax\ifhmode\unskip\space\fi \ZBLhref{#1}}
\providecommand{\ZBLhref}[1]{%
  \href{http://zbmath.org/?q=an:#1}{Zbl #1}}
\providecommand{\JfM}[1]{\relax\ifhmode\unskip\space\fi \JfMhref{#1}}
\providecommand{\JfMhref}[1]{%
  \href{http://zbmath.org/?q=an:#1}{JfM #1}}


\bibitem{MR0233275}
P.~Dembowski, \emph{Finite geometries}, Ergebnisse der {M}athematik und ihrer
  {G}renzgebiete ~44, Springer-Verlag, Berlin, 1968,
  \\
  \doi{10.1007/978-3-642-62012-6}. \MR{0233275}. \ZBL{0865.51004}.

\bibitem{MR0072144}
J.~A. Dieudonn{\'e}, \emph{La g\'eom\'etrie des groupes classiques}, Ergebnisse
  der {M}athematik und ihrer {G}renzgebiete ({N}.{F}.) ~5, Springer-Verlag,
  Berlin, 1955, \\
  \doi{10.1007/978-3-662-59144-4}. \MR{0072144}.
  \ZBL{0221.20056}.

\bibitem{MR1859189}
L.~C. Grove, \emph{Classical groups and geometric algebra}, Graduate Studies in
  Mathematics ~39, American Mathematical Society, Providence, RI,
  2002, \\
  \doi{10.1090/gsm/039}. \MR{1859189}. \ZBL{0990.20001}.

\bibitem{MR2795696}
T.~Grundh{\"o}fer, B.~Krinn, and M.~J. Stroppel, \emph{Non-existence of
  isomorphisms between certain unitals}, Des. Codes Cryptogr. \textbf{60}
(2011), no.~2, 197--201, \\
\doi{10.1007/s10623-010-9428-2}. \MR{2795696}.
  \ZBL{05909195}.

\bibitem{MR2926161}
M.~Gulde and M.~J. Stroppel, \emph{Stabilizers of subspaces under
  similitudes of the {K}lein quadric, and automorphisms of {H}eisenberg
  algebras}, Linear Algebra Appl. \textbf{437} (2012), no.~4, 1132--1161,
\\
\doi{10.1016/j.laa.2012.03.018}, \eprint{arXiv:1012.0502}. \MR{2926161}.
  \ZBL{06053093}.

\bibitem{MR514561}
S.~Helgason, \emph{Differential geometry, {L}ie groups, and symmetric spaces},
  Pure and Applied Mathematics ~80, Academic Press Inc., New York, 1978.
  \\
  \MR{514561 (80k:53081)}. \ZBL{0993.53002}.

\bibitem{MR0333959}
D.~R. Hughes and F.~C. Piper, \emph{Projective planes}, {G}raduate
  {T}exts in {M}athematics ~6, Springer-Verlag, New York, 1973. \MR{0333959}.
  \ZBL{0484.51011}.

\bibitem{MR0101253}
N.~Jacobson, \emph{Composition algebras and their automorphisms}, Rend. Circ.
  Mat. Palermo (2) \textbf{7} (1958), 55--80, \doi{10.1007/BF02854388}.
  \MR{0101253}. \ZBL{0083.02702}.

\bibitem{MR1009787}
N.~Jacobson, \emph{Basic algebra. {II}}, W. H. Freeman and Company, New York,
  2nd ed., 1989. \MR{1009787}. \ZBL{0694.16001}.

\bibitem{MR1187635}
B.~C. Kestenband, \emph{Generalizing a nonexistence theorem in finite
  projective planes}, Geom. Dedicata \textbf{44} (1992), no.~2, 123--126,
\\
\doi{10.1007/BF00182943}. \MR{1187635}. \ZBL{0772.51002}.

\bibitem{MR3535075}
N.~Knarr and M.~J. Stroppel, \emph{Heisenberg groups over composition
  algebras}, Beitr. Algebra Geom. \textbf{57} (2016), no.~3, 667--677,
\\
\doi{10.1007/s13366-015-0276-0}. \MR{3535075}. \ZBL{06619596}.

\bibitem{KramerStroppel-JoLT}
L.~Kramer and M.~J. Stroppel, \emph{Hodge operators and exceptional
  isomorphisms between unitary groups}, J. Lie Theory \textbf{33} (2023),
no.~1, 329--360,
\\
\eprint{arxiv:2208.11044},
  \urlprefix\url{https://www.heldermann.de/JLT/JLT33/jlt33.htm}.

\bibitem{MR0295934}
M.~E. O'Nan, \emph{Automorphisms of unitary block designs}, J. Algebra
  \textbf{20} (1972), 495--511, \doi{10.1016/0021-8693(72)90070-1}.
  \MR{0295934}. \ZBL{0241.05013}.

\bibitem{MR1384300}
H.~Salzmann, D.~Betten, T.~Grundh{\"o}fer, H.~H{\"a}hl, R.~L{\"o}wen, and
  M.~J. Stroppel, \emph{Compact projective planes}, de {G}ruyter Expositions in
  Mathematics ~21, Walter de Gruyter \& Co., Berlin, 1995,
  \doi{10.1515/9783110876833}. \MR{1384300}. \ZBL{0851.51003}.

\bibitem{MR1724629}
M.~J. Stroppel, \emph{Homogeneous symplectic maps and almost homogeneous
  {H}eisenberg groups}, Forum Math. \textbf{11} (1999), no.~6, 659--672,
\\
\doi{10.1515/form.1999.018}. \MR{1724629}. \ZBL{0928.22008}.

\bibitem{MR2241352}
M.~J. Stroppel and H.~{Van Maldeghem}, \emph{Automorphisms of unitals},
  Bull. Belg. Math. Soc. Simon Stevin \textbf{12} (2005), no.~5, 895--908,
  \\
  \urlprefix\url{http://projecteuclid.org/euclid.bbms/1136902624}.
  \\
  \MR{2241352}. \ZBL{1139.51002}.

\bibitem{MR1189139}
D.~E. Taylor, \emph{The geometry of the classical groups}, Sigma Series in Pure
  Mathematics ~9, Heldermann Verlag, Berlin, 1992. \MR{1189139}.
  \ZBL{0767.20001}.

\bibitem{MR0218489}
J.~Tits, \emph{Tabellen zu den einfachen {L}ie {G}ruppen und ihren
  {D}arstellungen}, Springer-Verlag, Berlin, 1967, \doi{10.1007/BFb0080324}.
  \MR{0218489}. \ZBL{0166.29703}.

\end{thebibliography}


\vfill\noindent
Markus J. Stroppel\\
LExMath,
Fakult\"at f\"ur Mathematik und Physik\\
Universit\"at Stuttgart\\
70550 Stuttgart,
Germany

\end{document}